\documentclass{article}
\usepackage{mathptmx}
\usepackage[letterpaper, margin=1in]{geometry}
\usepackage[utf8]{inputenc}
\usepackage[T1]{fontenc}
\usepackage{amsmath}
\usepackage{amsfonts}
\usepackage{amssymb}
\usepackage[version=4]{mhchem}
\usepackage{stmaryrd}
\usepackage{bbold}
\usepackage{hyperref}
\usepackage[vertfit]{breakurl} 

\usepackage[backend=biber,style=numeric, arxiv=format]{biblatex}
\hypersetup{colorlinks=true, linkcolor=blue, filecolor=magenta, urlcolor=cyan,}
\title{\large \textbf{Completeness and Well-Definability of a Provability Degree Measure in Sufficiently Powerful Formal Systems, and Finite-Time Effective Knowers }}

\author{Rohan Bahl\\
rohan.bahl@ieee.org}
\date{}

\begin{document}
\maketitle

\begin{abstract} 
We show that including degrees of a particular kind of provability in the search target for any theorem-prover in sufficiently powerful formal systems over finite-sized statements preserves well-definition and a sufficient consistency while establishing completeness. Moreover, the union of such degrees is isomorphic to such a system's $\aleph_{0}$ statements and permits the construction of a best-possible (up to a quadratic term) finite-time theorem prover, $\varphi^{\prime}$, while still subject to limitations in formal systems. These results, owing to the fact that $\varphi^{\prime}$ may arise through the behavior of any unbounded inductive computation, establish results on the behavior of a class of computational processes.
\end{abstract}

\section*{Introduction \& Motivation:}
Mathematical theorem-proving remains a preeminent field of the subject. The Church-Turing thesis posits the formalizability of all physical states and concepts in mathematical terms; theorem-proving enables one to draw conclusions from such data as well. The Hilbert Program aimed to establish the equivalence of true and provable statements, but Godel's Incompleteness Theorems \cite{Raatikainen} halted this ambition. Later, advances in computability theory, namely the Time Hierarchy Theorem \cite{Hennie} showed arbitrarily-long proof-times, and consequently proof-lengths.\\

A novel Hilbert's program was established a talk by Cai et al. \cite{Cai} , who questioned whether replacing the proof operator " $\vdash$ " with "provability degrees" would validate an "Ultimate Hilbert's Program". They concluded that, dropping the requirement of a uniform list of axiom systems, this would not work, i.e. there would still exist some unsolvable statements in that system according to their chosen measure of unprovability.\\

We consider a similar approach, establishing slightly adapted degrees of unprovability. In this, they are extended to tautologically encapsulate all finite-sized statements in the sufficiently-powerful formal system considered, forming a degree of  \textit{unknowability}. We prove two main theorems - that the degree of unknowability is finite for all finite statements through showing an equivalence between the $\aleph_{0}$ elements of each set, and that a slightly adapted proof procedure, based on Hutter's Universal search \cite{Hutter}, considering a target of such degrees of unknowability, always terminates in the consistent system, thus showing well-definedness of these operators. Our algorithm does not contradict Godel's First Incompleteness theorem and can take arbitrarily long to terminate on statements.\\

Semantics of this algorithm can arise through a wide variety of behaviors, including unbounded computation in Godel Machines \cite{Schmidhuber}, Solomonoff Induction \cite{Legg} and UCAI \cite{Katayama}. In summary, it establishes a bound on degrees of provability.

\section*{Preliminaries:}
Let us establish the structure and conventions for all sufficiently powerful formal systems in this document. Without loss of generality, they, and thus the argument, transfer to UTMs. Let F(L, D, A) be a Formal system $F$ with language $L$, deductive rules $D$ and axioms $A$. As always, $L(\Sigma, \in)$ denotes a language $L$ with alphabet $\Sigma$ and membership relation $\in$, while $\mathrm{D}\left(\left[p_{i} \in \mathrm{~F}(\mathrm{~L}) \Rightarrow q_{i} \in \mathrm{~F}(\mathrm{~L})\right.\right.$ $\forall \mathrm{i} \in[1, \mathrm{~m}]])$ is the deductive rule system with implicational chaining and $\mathrm{A}\left(\left[a_{i} \in \mathrm{~F}(\mathrm{~L}) \; \forall \mathrm{i} \in[1, n]\right]\right)$ are the initial axioms of $F$. We will consider, $S \in F(L)$ unless stated otherwise; $\mathbb{T}$ and $\mathbb{F}$ will designate True and False, respectively and $\{A, D, !, ~  \Rightarrow \} \in F(L)$.\\ \\
Define the proof operator as:\\
$\vdash: F$ $\times$ S $\times$ d $\in\{S,!S\} \rightarrow\left(d^{\prime} \in\{S,!S\},\; p[0 \ldots I-1]:\left(p[0]=S\; \& \; p[i] \Rightarrow p[i+1] \; \forall i \in[0, I-2] \; \& \; p[I-1]=d^{\prime}\right)\right.$ \\
It is equivalent to the standard proof operator, since $\vdash: F \times S= \; \vdash: F \times S \times \mathbb{T}$.\\
Then, $\operatorname{Con}(F) =\nexists \mathrm{S}:(\exists F \vdash \mathrm{~S} \; \& \; \exists \mathrm{~F} \vdash!\mathrm{S})$ becomes the consistency requirement.

\section*{Provability and the Theorem Prover:}
Introduce the following modified Hilbert-Bernays provability axioms \cite{Smith}. They are used without reference in all further derivations.

\begin{enumerate}
  \item $(F \vdash S)[d]=!S \Leftrightarrow(F \vdash!S)[d]=S$
  \item $\nexists(F \vdash S) \Leftrightarrow \nexists(F \vdash!S)$
  \item $\exists F \vdash(\exists(F \vdash S) \&(F \vdash S)[d]=\mathbb{T}) \Rightarrow \exists(F \vdash S)$
  \item $F \vdash(A=\mathbb{T} \Rightarrow B=\mathbb{T}) \Rightarrow((F \vdash A)[d]=\mathbb{T} \Rightarrow(\exists F \vdash A \Rightarrow \exists F \vdash B))$  \hspace{2em} (\^{}) 
\end{enumerate}

F will possess a canonical theorem-proving procedure adapted from Paulson \cite{Paulson}\\
\fbox{\parbox{\textwidth}{
Maintain the following lists for query efficiency (ts and fs denote true and false statement lists, respectively) \\

Initially, ts:= D $\cup \mathrm{A} \cup\{\mathbb{T}, \mathbb{F}\}$

\hspace{2em} ts: $=\{S: F \vdash S$ already $\} \cup$ fs 

\hspace{2em} fs: $=\{S: F \vdash!S$ already $\} \cup$ ts \\
Routine $\mathrm{F} \vdash \mathrm{S}$ :=\\
;loop:

\hspace{2em}; $\forall \mathrm{g} \in \mathrm{ts}:$

\hspace{4em};;;deductions=[g + $\left.a_{i}(\mathrm{~g}) \forall a_{i} \in \mathrm{~F}(\mathrm{~A})\right] / /$ Only the $a_{i}(\mathrm{~g})$ are shown in statement bank

\hspace{4em};;;ts += deductions

\hspace{4em};;;fs += !deductions

\hspace{2em};;check = prepare(S)

\hspace{2em};;if check $\in \mathrm{t} \in \mathrm{ts}$ ?:\(\hookleftarrow\)  $($check , ts[check]) // Can be any ts[check] which is consistent 

\hspace{2em};;if check $\in \mathrm{f} \in \mathrm{fs}$ ?:\(\hookleftarrow\)  (check , fs[check]) // Can be any fs[check] which is consistent 

;\(\hookleftarrow\)  // Should never be reached
}}
\fbox{\parbox{\textwidth}{%
Subroutine prepare(S):=; \(\hookleftarrow\) S
}}\\

The theorem-prover's runtime admits an exponential speedup via representing deductions as a digraph with each edge starting at a node and ending at the set of all nodes which are axioms (elements of $F(A)$) applied to the first node (two nodes are identical if they are syntactically equivalent). A proof of $F$ then consists of the subgraph with all nodes in the subgraph being proved from $F(A)$. The construction of  $\varphi^{\prime}$ then incorporates this method to achieve a best-possible runtime.

\section*{Limitations:}
Apart from Godel's Incompleteness Theorems, there must exist $S$ which does not have finite proof. Assume the opposite, namely $\mathbb{N} \ni|F \vdash S| \; \forall \; S$. Then applying the theorem prover would always terminate. But $\operatorname{Con}(\mathrm{F}) \in\{\mathrm{S}: \operatorname{Con}(\mathrm{F}) \Rightarrow \nexists \mathrm{F} \vdash \mathrm{S}\}$ \cite{Raatikainen}, so this is a contradiction. This leads to the main result.

\section*{Main result:}
\subsubsection*{\underline{Theorem M}:}
There exists an algorithm, $\varphi(\mathrm{F}, \mathrm{S})$ mapping $S$ to $f(F, S)$, where $f$ is defined at the end of the following section. Furthermore, there exists $\varphi^{\prime}: \mathrm{F} \times \mathrm{S} \rightarrow \varphi(\mathrm{F}, \mathrm{S})$ such that $\left|\# \varphi^{\prime}(\mathrm{F}, \mathrm{S})\right|=\mathrm{O}\left(\mid \# \varphi_{O}(\mathrm{~F}\right.$, $\left.S)\left.\right|^{2}\right) \forall S$, where \# denotes the computation trace and $\varphi_{O}$ is equivalent to $\varphi$, denoting the fastest algorithm over all the latter's arguments. .

\subsubsection*{Proof:}
We will prove \underline{Theorem M} in two parts. Firstly, we establish the existence of a candidate algorithm by \underline{Theorem 1} (which follows from \underline{Theorem 2} and \underline{Theorem 3}). Then we apply \underline{Theorem 4} to construct $\varphi^{\prime}$ and establish the lower bound.

\subsection*{Construction of the unknowability and unprovability operators:}
We motivate the construction of the unprovability and unknowability operators by the fact that some statements can be proven to have no proof, assuming Con(F). This also appears in \cite{Cai} and resembles \cite{user21820}.
\subsection*{{[ Degrees of ]} unprovability:}
Let the unprovability operator \cite{Verbrugge} be:

$$
\square: F \times S \rightarrow C o n(F) \Rightarrow \nexists F:(F(A)-\{\operatorname{Con}(F)\}) \vdash S
$$

More specifically, $\square$ is equivalent to finding whether $S$ is computable, or degree of unprovability as an ordinal, which can be obtained through the combination of operators below.\\\\
{Define the  \textit{degree of unprovability} as:}

$$
\square^{n \in F(L)}(F, S)=\square(F, \square(F, \ldots \square(F, S) \ldots) \text { applied } n \text { times }
$$

By \cite{Cai}, it is known that $\forall n \in \mathbb{N} \; \exists \; F, S: F \vdash \square^{n}(F, S)$ and by \cite{user21820}, all admissible ordinals in $F$ are contained in $\left\{n: \exists F, S: F \vdash \square^{n}(F, S)\right\}$. The argument may be extended to show that $\square^{m}\left(F, n: \mid \square^{n}\right.$ $(F, S) \mid \in \mathbb{N}$ ) exhibits the same behavior.
\subsection*{[Degrees of] unknowability:}
Define the\textit{ degree of reachability} as:

$$
\textnormal{o}: F \times S \rightarrow \min \left\{n \in F(L):\left|\square^{n}(F, S)\right| \in \mathbb{N}\right\}
$$ \\ 
Where $n \in F(L)$ follows from \cite{Cai}\\ \\
Define the \textit{degree of unknowability} as repeated application of the degree of reachability:

$$
o^{n \in \mathbb{N}}(\mathrm{~F}, \mathrm{~S})=\mathrm{o}(\mathrm{o}(\ldots \mathrm{o}(\mathrm{~F}, \mathrm{~S}) \ldots))
$$

$\mathbb{N} \ni n \; \forall \; S$ will be shown in \underline{Theorem 2}. Then let $f(F, S)=o^{\mathbb{N} \ni n}(F, S)$\\
From hereon, we will refer to n with the above meaning, unless indicated otherwise.

\section*{Establishing the existence of $\varphi$ :}
\underline{Theorem 1}: There is an algorithm $\varphi$, such that $\left.\varphi: \mathrm{F} \times \mathrm{S} \rightarrow\left(\mathrm{n}:\left|o^{n}(\mathrm{~F}, \mathrm{~S})\right| \in \mathbb{N}\right)\right)$\\ \\
\underline{Proof:}\\
Let $\varphi \; = \;\vdash$ (which is interpreted as the Theorem-Proving procedure) and prepare(S) : = f(F, S). \underline{Theorem 1} then follows from \underline{Theorem 2}, which claims $\mathbb{N} \ni|\varphi(F, S)| \; \forall \;(F, S)$ and \underline{Theorem 3}, which claims $\mathbb{N} \ni|\# \varphi(F, S)| \forall(F, S)$.\\ \\
\underline{Theorem 2}: The degree of unknowability is always finite for finite sentences,

$$
\exists k_{s} \in \mathbb{N}: o^{k_{s}}(\mathrm{~F}, \mathrm{~S})=0 \; \forall \mathrm{~S}
$$ \\ \\
\underline{Proof:} \\ \\
Prove a contradiction to the contrapositive. Let $\mathrm{O}(\mathrm{F}, \mathrm{S})=\left[o^{i}(\mathrm{~F}, \mathrm{~S}) \; \forall \mathrm{i} \leq \mathrm{n} \in \mathbb{N}\right]$. Clearly there is a bijection between $\{\mathrm{S}:|\mathrm{S}| \in \mathbb{N}\}=\aleph_{0}=\sum_{i \leq n}\left|o^{i}(\mathrm{~F}, \mathrm{~S})\right|$. But by applying the reasoning on this argument with $|O(F, S)| \notin \mathbb{N}$, then $O(F, S)=\left[o^{1}(F, S),|O(F, S)| \notin \mathbb{N}\right]$, contradiction $\blacksquare$ 

The reasoning can also show that $k_{s}=2$; this does not impact the algorithm.\\\\
\underline{Theorem 3}: The algorithm will always halt for every sentence\\
Let $\operatorname{Crit}(\mathrm{S})\; =|\# \varphi(\mathrm{~F}, \mathrm{~S})| \in \mathbb{N}$. Then this is equivalent to showing Crit := (Crit(S) $\forall \mathrm{S})$\\\\
\underline{Proof.}\\
For Sake of Contradiction, let counter(F) = \{S: !Crit(S)\}. |counter(F)| $=0 \Rightarrow$\underline{ Theorem 3}, so we assume $|\operatorname{counter}(F)| \neq 0$. We show that $\exists G(F):|\operatorname{counter}(F)| \neq 0 \Rightarrow G(F) \in$ counter(F).\\

Let $H(F):=($ first := min\{S: !Crit(S))\} if $\exists S:!$ Crit(S)) else $y)$, where $y \notin\{\mathbb{T}, \mathbb{F}\}$, and $G(F)=(H(F)=\mathbb{T})$. Consider $\varphi(F, G(F))$. Here, $\nexists S: \; !\operatorname{Crit}(S) \Rightarrow|\# \varphi(F, G(F))| \in \mathbb{N}$, but $\exists S: \operatorname{Crit}(S) \Rightarrow H(F) \in\{\mathbb{T}, \mathbb{F}\}$. Nonetheless, determining whether $G(F) \Rightarrow f(F, G(F)) \notin \mathbb{N}$. This is due to the well-ordering of degrees of unknowability; if the algorithm cannot halt on that statement, it is in power equivalent at least to a degree of unknowability for which the algorithm cannot halt (a tautology). But then $!\operatorname{Crit}\left(\mathrm{G}(\mathrm{F})\right.$), since the degree of unknowability, $\mathrm{O}\left(\mathrm{F}, \mathrm{G}(\mathrm{F})\right.$), is $k_{G(F)}=\; \mid[1, \mathrm{G}(\mathrm{F})$ : !Crit(F, $\left.\mathrm{G}(\mathrm{F}))\right]|=|[1$, $\min \{G(F): \; \mid \# \varphi(F, G(F) \mid=\max \{\mid \# \varphi(F, S \mid \forall S\}\}[G(F)][$ check $]] \mid=2$.\\

But, since $G(F) \in$ counter$(F) \; \& \;|\# \varphi(F, G(F))| \in \mathbb{N}$, this is a contradiction, and $|\operatorname{counter}(F)|=0$ . $\therefore \nexists \mathrm{S}:|\# \varphi(\mathrm{~F}, \mathrm{~S})| \notin \mathbb{N}$ $\blacksquare$ 

\section*{Extending $\varphi$ to an optimal $\varphi^{\prime}$ :}

We now construct the algorithm taking inspiration from \cite{Hutter}.\\\\
Define the dovetailing operator, interleaving a series of computations with each other, as:\\
dovetail: $\left(p_{1} \times \ldots \times p_{m}\right) \times \mathbb{N} \ni \mathrm{c} \leq \mathrm{m} \rightarrow \# p_{1}[\hookleftarrow] \times \ldots \times \# p_{m}[\hookleftarrow] $\\

Each step is performed on the specified subroutine of the pointer, as the pointer then moves to the next subroutine skipping all subroutines which have halted. The entire procedure halts when a target number of programs have halted.\\

Let $\varphi^{\prime}\left(F\right.$, S) := dovetail $\left(\varphi(\mathrm{F}, \mathrm{S}), *\left[\right.\right.$\!proof of $\varphi^{\prime \prime}(\mathrm{F}, \mathrm{S}) \; \& \; \operatorname{Con}($ proof\!$\left.\left.\left.) \forall \varphi^{\prime \prime} \in \mathrm{F}(\mathrm{L})\right\}\right], 1\right)$\\

Here, *[] represents the execution of each element of its argument in a new round of execution. $\varphi$ ' must alternate in execution between $\varphi(F, S)$ and the collective of the other arguments. By abuse of notation, we now consider Con(proof) to determine the logical consistency with the declared proof of knowability as its target. It is checked after execution termination of the respective $\varphi$ ". \\

Let $\varphi_{T}$ be the first subprogram to halt when executing $\varphi^{\prime}$ and providing $\varphi^{\prime}(F, S)=R$. Let $|R|$ be the statement's proof complexity (or equivalently size) produced, and let enc $\left(\varphi_{T}\right)$ be the subroutine's Godel Numbering. Then we have: \\\\
\underline{Theorem 4:}\\
$\left|\# \varphi^{\prime}(F, S)\right| \leq \min \left\{2|\# \varphi(F, S)|,\left(1+\varphi_{C}\right)^{2}\right\}$ where $\varphi_{C}=\operatorname{enc}\left(\varphi_{T}\right)+e$, and $\mathrm{e}=\left|\# \varphi_{T}\right|+|\# \operatorname{Con}(\mathrm{R})|$ denotes the shortest time to prove and verify $S$ (including to a quadratic factor) \\\\
\underline{Proof.}
We compute $\left|\# \varphi^{\prime}(\mathrm{F}, \mathrm{S})\right|$. We have $\mathrm{e}=\left|\# \varphi_{T}\right|+|\# \mathrm{Con}(\mathrm{R})|$ denote the total time taken to execute and verify consistency of $\varphi_{T}$ and its proof system. Then $\varphi^{\prime}$ has begun to execute another program, $\varphi_{C}$ : enc $\left(\varphi_{C}\right)=\operatorname{enc}\left(\varphi_{T}\right)+e$. Since it has not executed more than enc \!$\left(\varphi_{C}\right)$ - i steps for subprogram i , then the total number of computational steps are bounded by $\sum_{i=1}^{\varphi_{C}}\left(\varphi_{C}-\mathrm{i}\right)=\frac{\varphi_{c}\left(\varphi_{c}+1\right)}{2}$. Since $\varphi$ is executed each time whenever another subprogram in $\varphi^{\prime}$ is:

$$
\left|\# \varphi^{\prime}(\mathrm{F}, \mathrm{~S})\right| \leq \min \left\{|\# \varphi(\mathrm{~F}, \mathrm{~S})|, \frac{\varphi_{c}\left(\varphi_{c}+1\right)}{2}\right\}+\min \left\{|\# \varphi(\mathrm{~F}, \mathrm{~S})|, \frac{\varphi_{c}\left(\varphi_{c}+1\right)}{2}\right\}
$$ \\
where each of the identical terms' contributions are from $\varphi$ and the other programs, respectively.$\blacksquare$  \\\\
Axiom (\^{}) becomes:\\

$\mathrm{F} \vdash(A=\mathbb{T} \Rightarrow \mathrm{B}=\mathbb{T}) \Rightarrow\left((F \vdash A)[d]=\mathbb{T} \Rightarrow(\exists F \vdash A \Rightarrow \exists F \vdash B) \; \& \; \varphi^{\prime}(A)[t s[A][-1]]=\varphi^{\prime}(B)[t s[B][-1]]\right)$\\\\
This concludes the construction of $\varphi^{\prime}$.

\sloppy
\raggedright
\printbibliography

\end{document}